\documentclass[11pt]{article}

\usepackage{amsmath, amsthm, amssymb}
\usepackage{enumerate}
\usepackage{graphicx}
\usepackage{float}
\usepackage{hyperref}
\usepackage[margin=1in]{geometry}
\usepackage{comment}
\usepackage{color}
\numberwithin{equation}{section}

\theoremstyle{plain}
\newtheorem{theorem}{Theorem}[section]
\newtheorem*{theorem*}{Theorem}
\newtheorem{proposition}[theorem]{Proposition}
\newtheorem{lemma}[theorem]{Lemma}

\newtheorem*{conjecture*}{Conjecture}
\theoremstyle{definition}

\newtheorem{remark}[theorem]{Remark}

\begin{document}

% average integral symble
\def\Xint#1{\mathchoice
   {\XXint\displaystyle\textstyle{#1}}%
   {\XXint\textstyle\scriptstyle{#1}}%
   {\XXint\scriptstyle\scriptscriptstyle{#1}}%
   {\XXint\scriptscriptstyle\scriptscriptstyle{#1}}%
   \!\int}
\def\XXint#1#2#3{{\setbox0=\hbox{$#1{#2#3}{\int}$}
     \vcenter{\hbox{$#2#3$}}\kern-.5\wd0}}
\def\ddashint{\Xint=}
\def\dashint{\Xint-}

\title{A Liouville theorem for subcritical Lane-Emden system}

\author{Ze Cheng$^{a}$ \footnote{Partially supported by NSF DMS-1405175.},
 Genggeng Huang$^{b}$ \footnote{Partially supported by NSFC-11401376 and China Postdoctoral Science Foundation 2014M551391.},
 and Congming Li$^{b,a}$ \footnote{Partially supported by NSFC-11571233 and NSF DMS-1405175.}
 \\
$^a$ Department of Applied Mathematics, \\
University of Colorado Boulder, CO 80309, USA \\
$^b$ Department of Mathematics, INS and MOE-LSC, \\
Shanghai Jiao Tong University, Shanghai, China
}

\maketitle
\date{}

\begin{abstract}
In this paper, we present a necessary and sufficient condition to the Lane-Emden conjecture.
This condition is an energy type of integral estimate on solutions to subcritical Lane-Emden system. To approach the long standing and interesting conjecture, we believe that one plausible path is to refocus on establishing this energy type estimate. %In addition, a mixed type $W^{2,p}$-estimate is introduced to resolve some technical issue. %This Liouville type theorem can also be generalized.

%The existing proofs found by authors are effective under more or less stronger condition.

\end{abstract}

\section{Introduction}
This paper is devoted to the nonexistence of positive solution to the Lane-Emden system,
\begin{align}\label{laneEmdenSystem}
  \left\{\begin{array}{ll}
    %-\Delta u = c_1(x) v^p  & \text{in } \mathbb{R}^n, \\
    %-\Delta v = c_2(x) u^q  & \text{in } \mathbb{R}^n, \\
    -\Delta u = v^p  , \\ %& \text{in } \mathbb{R}^n, \\
    -\Delta v = u^q  , %& \text{in } \mathbb{R}^n, \\
	%\ \ u,v \geq 0 & \text{in } \mathbb{R}^n, \\
 	%u_1(x),u_2(x) \rightarrow 0   & \text{uniformly as } |x| \rightarrow +\infty.
        \end{array}
\right.
\quad \text{in } \mathbb{R}^n,
\end{align}
where $\ u,v \geq 0$, $0<p, q < +\infty$.
%, and $c_1(x), c_2(x)$ are double bounded, namely, there exists a constant $C>0$ such that
%\begin{align}\label{doubleBounded}
%	\frac{1}{C} \leq c_1(x),  c_2(x) \leq C, \ \forall x\in \mathbb{R}^n.
%\end{align}
%An important special case is Lane-Emden system, for which $c_1(x)=c_2(x)\equiv 1$.
%In this paper we consider the non-existence of solution to the system, and
The hyperbola \cite{Mitidieri93,Mitidieri96}
\begin{align*}
	\frac{1}{p+1}+\frac{1}{q+1}=\frac{n-2}{n}
\end{align*}
is called critical curve because it is known that on or above it, i.e.
\begin{align*}
	\frac{1}{p+1}+\frac{1}{q+1} \leq \frac{n-2}{n},
\end{align*}
which is called critical and  supercritical respectively,
the system \eqref{laneEmdenSystem} admits (radial) non-trivial solutions, cf. Serrin and Zou \cite{SZ98}, Liu, Guo and Zhang \cite{LGZ06} and Li \cite{Li13}.
However, for subcritical cases, i.e. $(p,q)$ satisfying,
\begin{align}\label{subcriticalRegion}
\frac{1}{p+1}+\frac{1}{q+1}>\frac{n-2}{n},
\end{align}
people guess that the following statement holds and call it the Lane-Emden conjecture: 
\begin{conjecture*}
 $u=v\equiv0$ is the unique nonnegative solution for system \eqref{laneEmdenSystem}.
\end{conjecture*}

The full Lane-Emden conjecture is still open.  Only partial results are known, and many researchers have made contribution in pushing the progress forward. We shall briefly present some important recent developments of the Lane-Emden conjecture.

%For $n\leq 4$, the conjecture is completely solved. For higher dimensions, people can prove the conjecture under some additional restrictions on $p,q$.
Denote the scaling exponents of system \eqref{laneEmdenSystem} by
\begin{align}\label{scalingComponent}
\alpha = \dfrac{2(p+1)}{pq-1}, \ \ \beta = \dfrac{2(q+1)}{pq-1},\quad \text{for }pq>1.
\end{align}
Then subcritical condition \eqref{subcriticalRegion} is equivalent to
\begin{align}\label{subcriticalRegion2}
	\alpha + \beta > n-2,\quad \text{for } pq>1.
\end{align}
 For $p,q$ in the following region
\begin{align}
	pq\leq 1, \ \text{or } pq>1 \ \text{and } \max\{\alpha,\beta\}\geq n-2,
\end{align}
\eqref{laneEmdenSystem} admits no positive entire \emph{supersolution}, cf. Serrin and Zou \cite{SZ96}. This implies the conjecture for $n=1,2$.
Also, the conjecture is true for
\begin{align}\label{BMRegion}
	\min\{\alpha,\beta\}\geq \frac{n-2}{2}, \ \text{with } (\alpha,\beta) \neq (\frac{n-2}{2},\frac{n-2}{2}),
\end{align}
cf. Busca and Man\'{a}sevich \cite{BM02}. Note that \eqref{BMRegion} covers the case that both $(p,q)$ are subcritical, i.e. $\max\{p,q\}\leq \frac{n+2}{n-2}$, with $(p,q)\neq (\frac{n+2}{n-2},\frac{n+2}{n-2})$, which is treated earlier, cf. de Figueiredo and Felmer \cite{FF94} and Reichel and Zou \cite{RZ00}.
Also, Mitidieri \cite{Mitidieri96} has proved that the system admits no radial positive solution.
Chen and Li \cite{CL09a} have proved that any solution with finite energy must be radial, therefore combined with Mitidieri \cite{Mitidieri96}, no finite-energy non-trivial solution exists.

For $n=3$,  the conjecture is solved by two papers. First, Serrin and Zou \cite{SZ96} proved that there is no positive solution with polynomial growth at infinity. Then Pol\'{a}\v{c}ik, Quittner and Souplet \cite{PQS07} removed the growth condition. In fact, they proved that no bounded positive solution implies no positive solution. This result has two important consequences. One is that combining with Serrin and Zou's result, one can prove the conjecture for $n=3$. The other is that proving the Lane-Emden conjecture is equivalent to proving nonexistence of bounded positive solution.
Thus, we always assume that $(u,v)$ are bounded in this paper.

For $n=4$, the conjecture is recently solved by Souplet \cite{Souplet09}. In \cite{SZ96}, Serrin and Zou used the integral estimates to derive the nonexistence results. Souplet further developed the approach of integral estimates and solved the conjecture for $n=4$ along the case $n=3$. In higher dimensions, this approach provides a new subregion where the conjecture holds, but the problem of full range in high dimensional space still seems stubborn.
Souplet has proved that if
\begin{align}\label{soupletRegion}
	\max\{\alpha,\beta\}>n-3,
\end{align}
then \eqref{laneEmdenSystem} with $(p,q)$ satisfying \eqref{subcriticalRegion} has no positive solution. Notice that \eqref{soupletRegion} covers \eqref{subcriticalRegion} only when $n\leq 4$, and when $n\geq 5$ \eqref{soupletRegion} covers a subregion of \eqref{subcriticalRegion}.

%Recently, in \cite{LeiLi2013} Lei and Li have proved that, for system \eqref{doubleBoundedSystem} to admit positive entire solution for some double bounded coefficients $c_1(x), c_2(x)$ if and only if $\max\{\alpha,\beta\}\leq n-2$.

%A graph here

%Notice that in the subcritical region, precisely, $\{\max\{\alpha,\beta\}\leq n-2\}\cap\{\alpha+\beta>n-2\}$ there exists solution to \eqref{doubleBoundedSystem} for some double bounded coefficients $c_1(x), c_2(x)$.

%Here is an example...

%Notice that the coefficient is radially decreasing. This indicates that the nature of the coefficient matters to the existence of positive entire solutions. In other word, to prove Lane-Emden conjecture, one needs to exploit the fact that the coefficients are constant, 1.

%1. solution from lower dimension to higher dimension.

%2. is it essentially the monotonicity or 1 that matters to the nonexistence?

The approach developed by Souplet in \cite{Souplet09} is  also effective on non-existence of positive solution to Hardy-H\'{e}non type equations and systems (cf. \cite{Fazly2014, FG2014, Phan12, PS12}):
\begin{equation*}
\begin{cases}-\Delta u=|x|^a v^p,\\
-\Delta v=|x|^b u^q,
\end{cases}\quad \text{ in }\mathbb R^n.
\end{equation*}
This approach can also be applied to more general elliptic systems, for further details, we refer to \cite{Souplet12} and \cite{QS12}. Moreover, a natural extension and application of this tool is the high order Lane-Emden system which was done by Arthur, Yan and Zhao \cite{AYZ14}.

In this paper, we point out that the key to the Lane-Emden conjecture is obtaining a certain type of energy estimate. This estimate is in fact a necessary and sufficient condition to the conjecture. Connecting the estimate and the conjecture is a laborious work and needs to combine many types of estimates. We believe that with the result here people can refocus on proving the crucial estimate and thus solve the conjecture.

\begin{theorem}\label{liouvilleThm}
Let $n\geq 3$ and $(u,v)$ be a non-negative bounded solution to \eqref{laneEmdenSystem}. %where $c_1(x), c_2(x)$ are radially non-decreasing, i.e.
%\begin{align*}
	%x\cdot \nabla c_1(x) \geq 0, \  x\cdot \nabla c_2(x) \geq 0, \ \forall x\in \mathbb{R}^n.
%\end{align*}
%Let $r=\frac{s(q+1)}{p+1}$, and if there exists an $s>0$ satisfying $n-s\beta=n-r\alpha<1$ such that
Assume %additionally $p\geq q$ and
there exists an $s>0$ satisfying $n-s\beta<1$ such that
\begin{align}
    \int_{B_R} v^s \leq CR^{n-s\beta},  \label{vEstimateBetter}
    %\int_{B_R} u^r\leq CR^{n-r\alpha}, \label{uEstimate}
\end{align}
then $u, v \equiv 0$ provided $0<p,q<+\infty$ and  $\frac 1{p+1}+\frac 1{q+1}>1-\frac 2n$.
\end{theorem}

\begin{remark}\label{remark2}
\begin{enumerate}
\item[(a)] Energy estimate \eqref{vEstimateBetter} is a necessary condition to the Lane-Emden conjecture. One just needs to notice that when $u,v\equiv 0$, \eqref{vEstimateBetter} is obviously satisfied. The key to the proof of Theorem \ref{liouvilleThm} is to show \eqref{vEstimateBetter} is sufficient.

\item[(b)] If $p\geq q$, the assumption on $v$ is weaker than the corresponding assumption on $u$ due to a comparison principle between $u$ and $v$ (i.e. Lemma \ref{comparisonPrinciple}).

In other words, if $p\geq q$, and we assume for some $r>0$, such that $n-r\alpha <1$,
\begin{align}
	\int_{B_R} u^r \leq CR^{n-r\alpha}.  \label{uEstimateBetter}
\end{align}
Then \eqref{uEstimateBetter} implies \eqref{vEstimateBetter} by Lemma \ref{comparisonPrinciple}.

%As in our proof, we shall assume $p\geq q$ and prove with this weaker form of energy estimate \eqref{vEstimateBetter}.

\item[(c)] By taking $s=p$ Theorem \ref{liouvilleThm} recovers the result in \cite{Souplet09}. 

%\item[(a)] The proof of Theorem \ref{liouvilleThm} relies on a `proportional' structure of the estimate, in fact it is an open problem if we can remove this structure (see Open Problem 2).
\item[(d)] A technical issue is that the standard $W^{2,p}$-estimate used in \cite{Souplet09} is not enough to establish Theorem \ref{liouvilleThm} (see the footnote of Proposition \ref{estimateOnBR}). To overcome this difficulty, a mixed type $W^{2,p}$-estimate is introduced in Lemma \ref{w2pEstimate}.
\end{enumerate}
\end{remark}

\begin{remark}
\begin{enumerate}
\item[(a)] It is worthy to point out an interesting role that the coefficient ``1" of the nonlinear source term plays in the Lane-Emden system. Consider the following system
\begin{equation}\label{laneEmdenSystem2}
\begin{cases}
-\Delta u=c_1(x)v^p,\\
-\Delta v=c_2(x)u^q,
\end{cases}\quad\text{in}\quad \mathbb R^n,
\end{equation}
where $0<a\leq c_1(x),c_2(x)\leq b<\infty$ and $x\cdot\nabla c_1(x),x\cdot\nabla c_2(x)\geq 0$ for some positive constants $a,b>0$. We can also have the following Rellich-Poho\v{z}aev type identity for some constants $d_1,d_2$ such that $d_1+d_2=n-2$,
\begin{equation}
\label{pohozaevId2}
\begin{split}
&\int_{B_R}(\frac{nc_1}{p+1}-d_1c_1+\frac{x\cdot \nabla c_1(x)}{p+1})v^{p+1}+(\frac{nc_2}{q+1}-d_2c_2+\frac{x\cdot \nabla c_2(x)}{q+1})u^{q+1}dx \\
		&= R^n \int_{\mathbb S^{n-1}} \frac{c_1(R)v^{p+1}(R)}{p+1}+\frac{c_2(R)u^{q+1}(R)}{q+1}d\sigma  \\
		& +R^{n-1} \int_{\mathbb S^{n-1}} d_1v'u+d_2u'v d\sigma +
		R^n\int_{\mathbb S^{n-1}}(v'u'-R^{-2}\nabla_{\theta}u\cdot\nabla_{\theta}v)d\sigma.	
\end{split}
\end{equation}
By the constrains on $c_1(x),c_2(x)$, we can have the left terms (LT) in \eqref{pohozaevId2} as
\begin{equation}
LT\geq \delta_0\int_{B_R}v^{p+1}+u^{q+1}dx, \quad \text{for some}\quad \delta_0>0.
\end{equation}
The argument in \cite{Souplet09} is also valid for this case, and we still can prove nonexistence for $n\leq 4$ and for $\max(\alpha,\beta)>n-3, n\geq 5$. 

On the other hand, for $c_1(x),c_2(x)$ such that $x\cdot \nabla c_1(x), x\cdot \nabla c_2(x)<0$, there exist non-zero solutions of \eqref{laneEmdenSystem2} in some subcritical cases (see Lei and Li \cite{LL13} for detail). %It indicates that the monotonicity of $c_1(x),c_2(x)$ in the radial direction plays an important role in the nonexistence of positive solution.

\item[(b)] Theorem \ref{liouvilleThm} is still true if we consider \eqref{laneEmdenSystem2} with $0<a\leq c_1(x),c_2(x)\leq b<\infty$ and $x\cdot \nabla c_1(x), x\cdot \nabla c_2(x)\geq 0$. And the proof is highly similar to the case $c_1=c_2=1$. So in this paper, we only prove for $c_1=c_2=1$.
\end{enumerate}
\end{remark}

The complete solution of the Lane-Emden conjecture may be a longstanding work. Hence, it will be interesting to consider the Lane-Emden conjecture under some conditions weaker than \eqref{vEstimateBetter}.
\\
\textbf{Open problem 1.} Can we prove the Lane-Emden conjecture under the following pointwise asymptotic:
\begin{equation*}
|v(x)|\leq C|x|^{-\gamma},\quad \text{for some}\quad 0<\gamma<\beta.
\end{equation*}
\\
\textbf{Open problem 2.} Can we prove the Lane-Emden conjecture under the following integral asymptotic:
\begin{equation*}
\int_{B_R}v^s\leq CR^\delta, \quad \text{for some} \quad s>0,\quad 0<\delta<1.
\end{equation*}
Clearly, if problem 2 is solved, problem 1 directly follows by choosing sufficiently large $s$. 

The paper is organized as follows. In Section 2, we provide a few preliminary results.  Some simplified proofs are given for the completeness and convenience of readers. One of the difficulty in the proof of Theorem \ref{liouvilleThm} was to control the embedding index, and we derived a varied form of $W^{2,p}$-estimate (see Lemma \ref{sphereEstimateW2p}) to solve this problem. In Section 3, we give the proof of Theorem \ref{liouvilleThm}. Our proof by classifying the argument into two cases hopefully can deliver the idea and the structure of the proof to readers in a clearer way.

\section{Preliminaries}
Throughout this paper, the standard Sobolev embedding on $\mathbb S^{n-1}$ is frequently used. Here we make some conventions about the notations.
Let $D$ denote the gradient with respect to standard metric on manifold. Let $n\geq 2$, $j\geq 1$ be integers and $1\leq z_1<\lambda\leq +\infty$, $z_2\neq (n-1)/j$. For $u=u(\theta)\in W^{j,z_1}(\mathbb S^{n-1})$, we have
\begin{align}\label{sobolevEmbedding}
	\|u\|_{L^{z_2}(\mathbb S^{n-1})} \leq C\left( \|D_{\theta}^j u\|_{L^{z_1}(\mathbb S^{n-1})} + \|u\|_{L^1(\mathbb S^{n-1})}\right) ,
\end{align}
where
\begin{align*}
\left\lbrace
\begin{array}{ll}
	\frac{1}{z_2} = \frac{1}{z_1} - \frac{j}{n-1}, \ &\text{if } z_1<(n-1)/j, \\
	z_2 = \infty, \ &\text{if } z_1>(n-1)/j,
\end{array}
\right.
\end{align*}
and $C=C(j,z_1,n)>0$. Although $C$ may be different from line to line, we always denote the universal constant by $C$. For simplicity, in what follows, for a function $f(r,\theta)$, we define
\begin{align}\label{sphereNorm}
\|f\|_{p}(r)=\|f(r,\cdot)\|_{L^p(\mathbb S^{n-1})},
\end{align}
if no risk of confusion arises. Also let $s,p,q$ be defined as in Theorem \ref{liouvilleThm} and
$$l=s/p, \quad k=\frac{p+1}{p},\quad m=\frac{q+1}q.$$
By Remark \ref{remark2} (b) and Lemma \ref{comparisonPrinciple}, throughout the paper, we always assume $p\geq q$.
At last, we set
$$F(R)=\int_{B_R}u^{q+1}dx.$$

\subsection{Basic Inequalities}
Let us start with a basic yet important fact. %which simply states a fact that $L^r$-norm on a layer of sphere $\partial B_{\tilde{R}}$ where $\tilde{R}\in [R,2R]$ can be estimated by $L^r$-norm on $B_{2R}$.
Considering $L^t$-norm on $B_{2R}$, we can write
\begin{align*}
	\|f\|_{L^t(B_{2R})}^t = \int_0^{2R} \|f(r)\|_{L^t(\mathbb S^{n-1})}^t r^{n-1} dr,
\end{align*}
then by a standard measurement argument (cf. \cite{SZ96}, \cite{Souplet09}) one can prove that:
\begin{lemma}\label{sphereEstimate}
Let $f_i\in L^{p_i}_{loc}(\mathbb{R}^n)$, and $i=1,\ldots,N$,  then for any $R>0$, there exists $\tilde{R}\in[R,2R]$  such that
\begin{align*}
	\|f_i\|_{L^{p_i}(\mathbb S^{n-1})}(\tilde{R})  \leq (N+1) R^{-\frac{n}{p_i}}\|f_i\|_{L^{p_i}(B_{2R})}, \ \text{for each } i=1,\ldots,N.
\end{align*}
\end{lemma}

The following lemma is a varied $W^{2,p}$-estimate which seems not to appear in any literature, so we give a simple proof.
\begin{lemma}\label{w2pEstimate}
Let $1<\gamma<+\infty$ and $R>0$. For $u\in W^{2,\gamma}(B_{2R})$, we have
\begin{align*}
	%\left( \int_{B_R} |D^2 u|^p \right) ^{\frac{1}{p}}\leq \left( \int_{B_{2R}} |\Delta u|^p\right) ^ {\frac{1}{p}}+ R^{\frac{n}{p}-(n+2)}\int_{B_{2R}} |u| \\
	\|D^2 u\|_{L^\gamma(B_R)} \leq C\left( \|\Delta u\|_{L^\gamma(B_{2R})} + R^{\frac{n}{\gamma}-(n+2)}\|u\|_{L^1(B_{2R})}\right)
\end{align*}
where $C=C(\gamma,n)>0$.
\end{lemma}
Proof. First we deal with functions in $C^2(B_2)\cap C^0(\overline{B_2})$. By standard elliptic $W^{2,p}$-estimate, we have
\begin{align}\label{w2pEstimateStd}
	\|D^2 u\|_{L^\gamma(B_1)} &\leq C( \|\Delta u\|_{L^\gamma(B_{\frac{3}{2}})}+\|u\|_{L^\gamma(B_{\frac{3}{2}})}).
\end{align}

By Lemma \ref{sphereEstimate}, $\exists \tilde{R}\in[\frac{7}{4},2]$ such that on $B_{\tilde{R}}$,
$u$ can be written as $u=w_1+w_2$, where respectively $w_1$ and $w_2$ are solutions to
\begin{align*}
\left\lbrace
\begin{array}{ll}
	\Delta w_1 = \Delta u, \ &\text{in } B_{\tilde{R}}, \\
	w_1=0, \ &\text{on } \partial B_{\tilde{R}},
\end{array}
\right.
\end{align*}
and
\begin{align*}
\left\lbrace
\begin{array}{ll}
	\Delta w_2 = 0, \ &\text{in } B_{\tilde{R}}, \\
	w_2=u, \ &\text{on } \partial B_{\tilde{R}},
\end{array}
\right.
\end{align*}
and additionally,
\begin{align}\label{sphereEstimateW2p}
	%\tilde R\int_{\partial B_{\tilde{R}}} |u| d\sigma\leq C \|u\|_{L^1(B_2)}\Rightarrow
	\int_{\partial B_{\tilde R}}ud\sigma \leq C\|u\|_{L^1(B_2)}.
\end{align}
By standard $W^{2,p}$-estimate with homogeneous boundary condition, we have
\begin{align*}
	\|w_1\|_{L^\gamma(B_{\frac{3}{2}})} \leq \| w_1\|_{W^{2,\gamma}(B_{\frac 32})}
			\leq C\|\Delta w_1\|_{L^\gamma(B_{\tilde{R}})}.
\end{align*}
Since $w_2$ can be solved explicitly by Poisson formula on $B_{\tilde{R}}$, we see that by \eqref{sphereEstimateW2p} for any $x\in B_{\frac{3}{2}}\subsetneq B_{\tilde{R}}$, $w_2(x)$ can be bounded pointwisely by
\begin{align*}
	|w_2(x)| \leq C \int_{\partial B_{\tilde{R}}} |u| \leq C\|u\|_{L^1(B_2)}.
\end{align*}
So,
\begin{align*}
	\|w_2\|_{L^\gamma(B_{\frac{3}{2}})} \leq C \|u\|_{L^1(B_2)}.
\end{align*}
Hence,
\begin{align*}
	\|u\|_{L^\gamma(B_{\frac{3}{2}})} &\leq \|w_1\|_{L^\gamma(B_{\frac{3}{2}})}+\|w_2\|_{L^\gamma(B_{\frac{3}{2}})} \\
			&\leq C(\|\Delta u\|_{L^\gamma(B_{\tilde{R}})}+\|u\|_{L^1(B_2)}).
\end{align*}
Therefore, \eqref{w2pEstimateStd} becomes
\begin{align*}
	\|D^2 u\|_{L^\gamma(B_1)} &\leq C( \|\Delta u\|_{L^\gamma(B_2)}+\|u\|_{L^1(B_2)}).
\end{align*}
Then the lemma follows from a dilation and approximation argument.
$\Box$

\begin{lemma}[Interpolation inequality on $B_R$]\label{interpolation}
Let $1\leq \gamma<+\infty$ and $R>0$. For $u\in W^{2,\gamma}(B_R)$, we have
\begin{align*}
	\|D_x u\|_{L^1(B_R)} \leq C\left( R^{n(1-\frac{1}{\gamma})+1} \|D_x^2 u\|_{L^\gamma(B_R)} + R^{-1} \|u\|_{L^1(B_R)}\right),
\end{align*}
where $C=C(\gamma,n)>0$.
\end{lemma}

\subsection{Poho\v{z}aev Identity, Comparison Principle and Energy Estimates}

For system \eqref{laneEmdenSystem} we have a Rellich-Poho\v{z}aev identity, which is the starting point of the proof of Theorem \ref{liouvilleThm}, %as in controlling $L^{q+1}$-norm of $u$ by quantities on $\mathbb{S}^{n-1}$.
\begin{lemma}\label{pohozaevId}
Let $d_1,d_2\geq0$ and $d_1+d_2=n-2$, then
\begin{align*}
		%&\int_{B_R}(\frac{nc_1}{p+1}-d_1c_1+\frac{x\cdot \nabla c_1(x)}{p+1})v^{p+1}+(\frac{nc_2}{q+1}-d_2c_2+\frac{x\cdot \nabla c_2(x)}{q+1})u^{q+1}dx \\
		%&= R^n \int_{\mathbb S^{n-1}} \frac{c_1(R)v^{p+1}(R)}{p+1}+\frac{c_2(R)u^{q+1}(R)}{q+1}d\sigma  \\
		%& +R^{n-1} \int_{\mathbb S^{n-1}} d_1v'u+d_2u'v d\sigma +
		%R^n\int_{\mathbb S^{n-1}}(v'u'-R^{-2}\nabla_{\theta}u\cdot\nabla_{\theta}v)d\sigma.		
		&\int_{B_R}(\frac{n}{p+1}-d_1)v^{p+1}+(\frac{n}{q+1}-d_2)u^{q+1}dx  \\
		&= R^n \int_{\mathbb S^{n-1}} \frac{v^{p+1}(R)}{p+1}+\frac{u^{q+1}(R)}{q+1}d\sigma
		 +R^{n-1} \int_{\mathbb S^{n-1}} d_1v'u+d_2u'v d\sigma +
		R^n\int_{\mathbb S^{n-1}}(v'u'-R^{-2}\nabla_{\theta}u\cdot\nabla_{\theta}v)d\sigma.
\end{align*}
\end{lemma}

\begin{lemma}[Comparison Principle]\label{comparisonPrinciple}
Let $p\geq q >0,pq>1$ and $(u,v)$ be a positive bounded solution of \eqref{laneEmdenSystem}. Then we have the following comparison principle,
\begin{align*}
	v^{p+1}(x)\leq \frac{p+1}{q+1} u^{q+1}(x), \  x\in \mathbb{R}^n.
\end{align*}
\end{lemma}
Proof. Let $l=(\frac{p+1}{q+1})^{\frac{1}{p+1}}$, $\sigma=\frac{q+1}{p+1}$. So $l^{p+1}\sigma=1$, and $\sigma\leq 1$. Denote
\begin{align*}
	\omega = v- lu^{\sigma}.
\end{align*}
We will show that $\omega \leq 0$.
\begin{align*}
	\Delta\omega &= \Delta v - l \nabla\cdot(\sigma u^{\sigma-1}\nabla u)  \\
				&= \Delta v - l\sigma(\sigma-1)|\nabla u|^2 -l\sigma u^{\sigma-1}\Delta u\\
				&\geq -u^q+l\sigma u^{\sigma-1}v^p \\
				&= u^{\sigma-1}((\frac{v}{l})^p - u^{q+1-\sigma}) \\
				&= u^{\sigma-1}((\frac{v}{l})^p - u^{\sigma p}).				
\end{align*}
So, $\Delta\omega>0$ if $w>0$.
Now, suppose $w>0$ for some $x\in\mathbb{R}^n$, and there are two cases:

Case 1: $\exists x_0\in \mathbb{R}^n$, such that $\omega(x_0)=\displaystyle \max_{\mathbb{R}^n} \omega(x)>0$, and $\Delta \omega(x_0)\leq 0$. However, when $w>0$, $\Delta\omega>0$, a contradiction.

Case 2: There exists a sequence $\{x_m\}$ with $|x_m|\rightarrow+\infty$, such that $\displaystyle \lim_{m\rightarrow+\infty} \omega(x_m) = \displaystyle \max_{\mathbb{R}^n} \omega(x) >c_0>0$ for some constant $c_0$.

Let $\omega_R(x)=\phi(\frac{x}{R})\omega(x)$, where $\phi(x)\in C_0^{\infty}(B_1)$ is a cutoff function and $\phi(x)\equiv1$ in $B_{\frac{1}{2}}$. Since $\omega_R(x)=0$ on $\partial B_R$, there exists an $x_R\in B_R$ such that $\omega_R(x_R)=\displaystyle \max_{B_R}\omega_R(x)$ and $\displaystyle \lim_{R\rightarrow+\infty} \omega(x_R) = \displaystyle \max_{\mathbb{R}^n} \omega(x) >0$. Also,
\begin{align*}
	0=\nabla \omega_R(x_R) = \phi(\frac{x_R}{R})\nabla\omega(x_R) + \frac{1}{R}\nabla\phi(\frac{x_R}{R})\omega(x_R).
\end{align*}
As $\phi(\frac{x_R}{R})\geq c_1>0$ for some constant $c_1$ (in fact, $\phi(\frac{x_R}{R})\rightarrow 1$) and $\omega(x_R)$ is bounded since $u,v$ are bounded in $\mathbb{R}^n$, we see that $\nabla\omega(x_R)\rightarrow 0$ as $R\rightarrow +\infty$.
So,
\begin{align*}
	0 &\geq \Delta\omega_R(x_R)=\frac{1}{R^2}\Delta\phi(\frac{x_R}{R})\omega(x_R)+\frac{2}{R}\nabla\phi(\frac{x_R}{R})\cdot\nabla\omega(x_R)+\phi(\frac{x_R}{R})\Delta\omega(x_R) \\
	\Rightarrow 0 &\geq \Delta\omega(x_R) + O(\frac{1}{R^2})
\end{align*}
Since $\omega(x_R)>c_0/2$ for sufficiently large $R$, $\Delta\omega(x_R)>c_2>0$ for some constant $c_2$, a contradiction.
$\Box$

\begin{remark}
For general Lane-Emden type system \eqref{laneEmdenSystem2}, we can choose
$$w=v-Clu^{\sigma},\quad \text{where}\quad C^{p+1}=\sup_{x\in\mathbb R^n}\frac{c_2(x)}{c_1(x)}.$$
By the same arguments, we can also get the desired comparison principle.

\end{remark}

Next we prove a group of energy estimates which are crucial to the entire argument in this paper. As Theorem \ref{liouvilleThm} points out, better energy estimates are the key to the Lane-Emden conjecture. Unfortunately, efforts have been made so far only provide the following inequalities, which are first obtained by Serrin and Zou \cite{SZ96} (1996). Here we give a simpler proof than the original one for the convenience of readers.

\begin{lemma}\label{energyEstimates}
Let $p,q>0$ with $pq>1$. For any positive solution $(u,v)$ of \eqref{laneEmdenSystem}
\begin{align}\label{uvEstimate1}
	\int_{B_R} u \leq C R^{n-\alpha}, \ \text{and } \int_{B_R} v \leq C R^{n-\beta},
\end{align}
\begin{align}\label{uvEsitmatePQ}
	\int_{B_R} u^q \leq C R^{n-q\alpha}, \ \text{and } \int_{B_R} v^p \leq C R^{n-p\beta}.
\end{align}
\end{lemma}

Proof. Without loss of generality, we can assume that $p\geq q$.

%Consider the following problem
%\begin{align*}
%\left\lbrace
%\begin{array}{ll}
%	-\Delta \phi = \lambda \phi \ &\text{in } B_R(0), \\
%	\phi = 0 \ &\text{on } \partial B_R(0).
%\end{array}
%\right.
%\end{align*}

Let $\phi\in C^{\infty}(B_R(0))$ be the first eigenfunction of $-\Delta$ in $B_R$ and $\lambda$ be the eigenvalue. By definition and rescaling, it is easy to see  that $\phi\mid_{\partial B_R} =0$ and $\lambda\sim \frac{1}{R^2}$. By normalizing, one gets $\phi\geq c_0>0$ on $B_{R/2}$ for some constant $c_0$ independent of $R$, $\phi(0)=\|\phi\|_{\infty}=1$. So, multiplying \eqref{laneEmdenSystem} by $\phi$ then integrating by parts on $B_R$ we have,
\begin{align*}
	\int_{B_R} \phi u^q = -\int_{B_R} \phi \Delta v &= \int_{\partial B_R} v \frac{\partial \phi}{\partial n} d\sigma + \lambda \int_{B_R} \phi v.
\end{align*}
By Hopf's Lemma we know that $\frac{\partial \phi}{\partial n}<0$ on $\partial B_R$, so
\begin{align}\label{uqBoundedByV}
	\int_{B_R} \phi u^q \leq \lambda \int_{B_R} \phi v.
\end{align}
Similarly, we have
\begin{align}\label{vpBoundedByU}
	\int_{B_R} \phi v^p \leq \lambda \int_{B_R} \phi u.
\end{align}
Applying Lemma \ref{comparisonPrinciple} to \eqref{uqBoundedByV}, we have
\begin{align*}
	\frac{1}{R^2} \int_{B_R} \phi v\geq C \int_{B_R} \phi v^{\frac{q(p+1)}{q+1}}.
\end{align*}
Notice that $\frac{q(p+1)}{q+1}>1$ as $pq>1$, so by H\"{o}lder inequality
\begin{align*}
	\int_{B_R} \phi v^{\frac{q(p+1)}{q+1}} &\geq (\int_{B_R} \phi v)^{\frac{q(p+1)}{q+1}}(\int_{B_R} \phi)^{-(\frac{q(p+1)}{q+1}-1)} \\
			&\geq C(\int_{B_R} \phi v)^{\frac{q(p+1)}{q+1}} R^{-n\frac{qp-1}{q+1}}.
\end{align*}
So,
\begin{align*}
	&\frac{1}{R^2} \int_{B_R} \phi v\geq C (\int_{B_R} \phi v)^{\frac{q(p+1)}{q+1}} R^{-n\frac{qp-1}{q+1}} \\
	\Rightarrow & \int_{B_R} \phi v \leq C R^{n-\beta}.
\end{align*}
Therefore, by \eqref{uqBoundedByV}
\begin{align*}
	\int_{B_R} \phi u^q \leq C R^{n-\beta-2} =CR^{n-q\alpha}.
\end{align*}
Now, \textbf{Case 1:} If $q>1$, then by H\"{o}lder inequality
\begin{align*}
	\int_{B_R} \phi u \leq (\int_{B_R} \phi u^q)^{\frac{1}{q}}(\int_{B_R} \phi )^{\frac{1}{q'}}
			\leq CR^{\frac{n}{q}-\alpha} R^{\frac{n}{q'}}
			= CR^{n-\alpha}, \quad \frac 1q+\frac 1{q'}=1.
\end{align*}
Mean while, by \eqref{vpBoundedByU}
\begin{align*}
	\int_{B_R} \phi v^p \leq CR^{n-\alpha-2} = CR^{n-p\beta}.
\end{align*}
This finishes the proof for Case 1.

\textbf{Case 2:} Assume that $q\leq 1$. To prove this trickier case, we begin with a lemma of energy-type estimate,

\begin{lemma}
If $\Delta u\leq 0$, then for $\gamma\in(0,1)$, $\eta\in C^\infty_0(\mathbb{R}^n)$,
\begin{align}\label{energy}
	\int_{\mathbb{R}^n} \frac{4}{\gamma^2}|D(u^{\frac{\gamma}{2}})|^2\eta^2 = \int \eta^2 |Du|^2 u^{\gamma-2}
				 \leq C \int |D\eta|^2u^{\gamma}.
\end{align}
\end{lemma}

Proof. Multiply $\eta^2 u^{\gamma-1}$ to $\Delta u\leq 0$ then integrate over the whole space. $\Box$

We rewrite \eqref{energy} as
\begin{equation}\label{energy2}
\int_{B_R}|D u|^2u^{\gamma-2}\leq \frac{C_\gamma}{R^2}\int_{B_{2R}}u^{\gamma}
\end{equation}
where $C_\gamma\rightarrow+\infty$ as $\gamma\rightarrow 1$.
From Poincar\'e's Inequality, we have
\begin{align}\label{embedding}
	|f|_{\frac{na}{n-a},\Omega_R} \leq C(n,a,\Omega)\left(|Df|_{a,\Omega_R} + |R|^{\frac{n-a}{a}}|f_{\Omega_R}|\right),
\end{align}
where
\begin{align*}
		f_{\Omega_R}=\dashint_{\Omega_R} f=\frac{1}{|\Omega_R|}\int_{\Omega_R} f,\quad \Omega_R=\{Rx|x\in\Omega\}.
\end{align*}
\par
Next we prove a variation of embedding inequality,
\begin{lemma}
For any $l\geq 1$,
\begin{align}\label{variationEmbedding}
	|f^l|_{\frac{an}{n-a},\Omega_R} \leq C(n,a,\Omega)\left( |D(f^l)|_{a,\Omega_R} + |R|^{\frac{n-a}{a}} |f_{\Omega_R}|^l\right)
\end{align}
\end{lemma}
Proof. By \eqref{embedding},
\begin{align*}
	|f^l|_{\frac{an}{n-a},\Omega_R} &\leq C(n,a,\Omega)\left( |D(f^l)|_{a,\Omega_R} + |R|^{\frac{n-a}{a}} |(f^l)_{\Omega_R}|\right)  \\
				&\leq  C(n,a,\Omega)\left( |D(f^l)|_{a,\Omega_R} + |R|^{\frac{n-a}{a}-n} \int_{\Omega_R} f^l dx \right) \\		
				&\leq  C(n,a,\Omega)\left\lbrace  |D(f^l)|_{a,\Omega_R} + |R|^{\frac{n-a}{a}-n} (\int_{\Omega_R} f dx)^{\theta l} (\int_{\Omega_R} f^{l\frac{an}{n-a}}dx)^{(1-\theta)l \frac{n-a}{lan}}\right\rbrace,\quad\theta=\frac{1-\frac{n-a}{na}}{l-\frac{n-a}{na}} \\
&\leq C(n,a,\Omega) |D(f^l)|_{a,\Omega_R}  + \frac 1 2 |f^l|_{\frac{an}{n-a},\Omega_R} + C(n,a,\Omega) |R|^{\frac{n-a}{a}} |f_{\Omega_R}|^l.	
\end{align*}
In getting the last inequality, we have used the Young's inequality. So we get \eqref{variationEmbedding}.
$\Box$

Let $l\geq 1$, $\theta\leq 2q<2$, $\gamma=l\theta<1$, $f=u^{\frac\theta2}$, $a=2$. Then
\begin{equation}
\begin{split}
|f^l|_{\frac{2n}{n-2},B_R}&\leq C\left(|Df^l|_{2,B_R}+R^{\frac{n-2}2}|f_{B_R}|^l\right)\\
&\leq C\left(|D(u^{\frac{l\theta}{2}})|_{2,B_R}+R^{\frac{n-2}{2}}|(u^{\frac\theta2})_{B_R}|^l\right)\\
&\leq \frac CR\left(\int_{B_{2R}}u^{l\theta}\right)^{\frac 12}+R^{\frac{n-2}2}\left(\dashint_{B_R}u^{\frac\theta2}\right)^l.
\end{split}
\end{equation}
The last term on the right can be estimate by H\"{o}lder and the fact that $\dashint_{B_R} u^q\leq CR^{-q\alpha}$ since $\frac{\theta}{2}<q$. This yields that
\begin{equation}
\int_{B_R}u^{\frac n{n-2}\theta l}\leq C\left(R^{-\frac{2n}{n-2}}\left(\int_{B_{2R}}u^{l\theta}\right)^{\frac n{n-2}}+R^{n-\frac n{n-2}l\theta\alpha}\right).
\end{equation}
This means if $\dashint_{B_R}u^{l\theta}\leq CR^{-l\theta\alpha}$, we have $\dashint_{B_R}u^{\frac n{n-2}l\theta}\leq CR^{-\frac{n}{n-2}l\theta\alpha}$ provided $l\theta<1$. By $\dashint_{B_R} u^q\leq CR^{-q\alpha}$, one gets
\begin{equation}
\dashint_{B_R}u^s\leq C(s)R^{-s\alpha},\quad \text{for}\quad s<\frac n{n-2}
\end{equation}
where $C(s)\rightarrow+\infty$ as $s\rightarrow \frac{n}{n-2}$.

By taking $s=1$, the above inequality immediately leads to $$\int_{B_R} u \leq C R^{n-\alpha}.$$
Since $pq>1$ and we assume that $p\geq q$, $q$ must be greater than 1, then by H\"{o}lder and \eqref{vpBoundedByU} we get
$$\int_{B_R} v^p \leq C R^{n-p\beta}.$$
This finishes the proof of Lemma \ref{energyEstimates}. $\Box$

\subsection{Key Estimates on $\mathbb{S}^{n-1}$}

Now that we have energy inequalities \eqref{uvEsitmatePQ}, in our assumption \eqref{vEstimateBetter} we can always assume $s\geq p$.
Since $l=\frac{s}{p}$, we have $l\geq 1$. The following estimates for quantities on sphere $\mathbb S^{n-1}$ are necessary to the proof.
\begin{proposition}\label{estimateOnSphere}
For $R\geq 1$, there exists $\tilde{R}\in[R,2R]$ such that for $l=\frac s p\geq 1$, $k=\frac{p+1}{p}$ and $m=\frac{q+1}{q}$, we have
\begin{align*}
	\|u\|_1(\tilde{R}) \leq CR^{-\alpha}, & \
	\|v\|_1(\tilde{R}) \leq CR^{-\beta},  \\
    \|D^2_x u\|_l(\tilde{R}) \leq CR^{-\frac{lp\beta}{l+\varepsilon}}, &\
       % \|D^2_x u\|_{1+\varepsilon}(\tilde{R}) &\leq CR^{-\frac{p\beta}{1+\varepsilon}}, \
        \|D^2_x v\|_{1+\varepsilon}(\tilde{R}) \leq CR^{-\frac{q\alpha}{1+\varepsilon}}, \\
        \|D_x u\|_1(\tilde{R}) \leq CR^{1-\frac{\alpha+2}{1+\varepsilon}},  &\
            \|D_x v\|_1(\tilde{R}) \leq CR^{1-\frac{\beta+2}{1+\varepsilon}}, \\
            \|D^2_x u\|_k(\tilde{R}) \leq  C(R^{-n}F(2R))^{\frac{1}{k}},&\
            \|D^2_x v\|_m(\tilde{R}) \leq  C(R^{-n}F(2R))^{\frac{1}{m}}.
\end{align*}
\end{proposition}

In view of Lemma \ref{sphereEstimate}, to prove Proposition \ref{estimateOnSphere}, we shall give the corresponding estimates on $B_{2R}$ first. We use the varied $W^{2,p}$-estimate (i.e. Lemma \ref{w2pEstimate}) to achieve this.
\begin{proposition}\label{estimateOnBR}
For $R\geq 1$, we have
\begin{align}\label{1stEstimate}
\left\{\begin{array}{ll}
	\|u\|_{L^1(B_R)} &\leq CR^{n-\beta}, \\
	\|v\|_{L^1(B_R)} &\leq CR^{n-\alpha},
	\end{array}
\right.
\end{align}

\begin{align}\label{2ndEstimate}
\left\{\begin{array}{ll}
	%\int_0^R \|D_x^2 u(r)\|_1 r^{n-1} dr &\leq CR^{n-p\beta}, \ R\geq 1, \\
	%\int_0^R \|D_x^2 v(r)\|_1 r^{n-1} dr &\leq CR^{n-q\alpha}, \ R\geq 1,
	\|D_x^2 u\|^{l+\varepsilon}_{L^{l+\varepsilon}(B_R)} \leq CR^{n-lp\beta}, \ \text{with } l=\frac s p\geq 1, \\
	%\|D_x^2 u\|^{1+\varepsilon}_{L^{1+\varepsilon}(B_R)} \leq CR^{n-p\beta}, \\
	\|D_x^2 v\|^{1+\varepsilon}_{L^{1+\varepsilon}(B_R)} \leq CR^{n-q\alpha},
	\end{array}
\right.
\end{align}

\begin{align}\label{3rdEstimate}
\left\{\begin{array}{ll}
	%R^{-1}\|D_x u\|_{L^1(B_R)} &\leq C\left( R^{n(1-\frac{1}{r})}\|D_x^2 u\|_{L^r(B_R)} + R^{-2} \|u\|_{L^1(B_R)}\right) , \ R\geq 1, \\
	%R^{-1}\|D_x v\|_{L^1(B_R)} &\leq C\left( R^{n(1-\frac{1}{t})}\|D_x^2 v\|_{L^t(B_R)} + R^{-2} \|v\|_{L^1(B_R)}\right), \ R\geq 1,
	\|D_x u\|_{L^1(B_R)} &\leq C R^{n+1-\frac{\alpha+2}{1+\varepsilon}},\\
	\|D_x v\|_{L^1(B_R)} &\leq C R^{n+1-\frac{\beta+2}{1+\varepsilon}},
	\end{array}
\right.
\end{align}
and let $k=\frac{p+1}{p}$, $m=\frac{q+1}{q}$,
\begin{align}\label{4thEstimate}
\left\{\begin{array}{ll}
	%\int_0^R \|D_x^2 u(r)\|_k^k r^{n-1} dr \leq CF(2R), \ R\geq 1,
	\|D_x^2 u\|_{L^k(B_R)}^k \leq CF(2R), \\
	%\int_0^R \|D_x^2 v(r)\|_m^m r^{n-1} dr \leq CF(2R), \ R\geq 1,
	\|D_x^2 v\|_{L^m(B_R)}^m \leq CF(2R).
\end{array}
	\right.
\end{align}
\end{proposition}
Proof.
Some frequently used facts include, $q\alpha=\beta+2$, $p\beta=\alpha+2$ and hence $n-kp\beta<0$ (due to \eqref{subcriticalRegion2}) and therefore $l<k$ (since $n-lp\beta\geq 0$).

Estimates \eqref{1stEstimate}  directly follow from \eqref{uvEstimate1} in Lemma \ref{energyEstimates}.

For the first estimate of \eqref{2ndEstimate}, after applying Lemma \ref{w2pEstimate}, the mixed type  $W^{2,p}$-estimate\footnote{Notice that with the standard $W^{2,p}$-estimate, we end up with a term of $\|u\|_{l+\epsilon}$ which cannot be estimated by any energy bound.}, we get
\begin{align*}
\|D_x^2 u\|^{l+\varepsilon}_{L^{l+\varepsilon}(B_R)} &\leq C\left( \|\Delta u\|^{l+\varepsilon}_{L^{l+\varepsilon}(B_{2R})} + R^{n-(l+\varepsilon)(n+2)}\|u\|^{l+\varepsilon}_{L^1(B_{2R})}\right).
\end{align*}
Then we use the assumed estimate \eqref{vEstimateBetter}  and Lemma \ref{energyEstimates} to get
\begin{align*}
	\|D_x^2 u\|^{l+\varepsilon}_{L^{l+\varepsilon}(B_R)} %&\leq C\left( \|\Delta u\|^{l+\varepsilon}_{L^{l+\varepsilon}(B_{2R})} + R^{n-(l+\varepsilon)(n+2)}\|u\|^{l+\varepsilon}_{L^1(B_{2R})}\right)  \\
				&\leq C\left( \int_{B_{2R}} v^{p(l+\varepsilon)}dx + R^{n-(l+\varepsilon)(n+2)}R^{(l+\varepsilon)(n-\alpha)}\right)  \\
				&\leq C\left( R^{n-pl\beta}+R^{n-(l+\varepsilon)(2+\alpha)}\right) \\
				&\leq C R^{n-pl\beta},
\end{align*}
where the last inequality is due to $\alpha+2=p\beta$. 
%Actually, in the argument above, $l$ can be replaced by any $\tilde{l}\in[1,l]$.
For the second estimate of \eqref{2ndEstimate},
\begin{align*}
	\|D_x^2 v\|^{1+\varepsilon}_{L^{1+\varepsilon}(B_R)} &\leq C\left(  \|\Delta v\|^{1+\varepsilon}_{L^{1+\varepsilon}(B_{2R})}  + R^{n-(1+\varepsilon)(n+2)}\|v\|^{1+\varepsilon}_{L^1(B_{2R})} \right) \\
			&\leq C\left( \int_{B_{2R}} u^{q(1+\varepsilon)}dx + R^{n-(1+\varepsilon)(n+2)} R^{(1+\varepsilon)(n-\beta)} \right) \\
			&\leq C\left( R^{n-q\alpha} + R^{n-(1+\varepsilon)(\beta+2)} \right) \\
			&\leq C R^{n-q\alpha}.
\end{align*}
%Similarly, we get estimate for $\|D_x^2 u\|^{1+\varepsilon}_{L^{1+\varepsilon}(B_R)}$.
For the first estimate of \eqref{3rdEstimate}, by Lemma \ref{interpolation}, % and by an argument similar to that we do to the 1st estimate in \eqref{2ndEstimate} just by replacing $l$ with $1+\varepsilon$ we see that
\begin{align*}
	\|D_x u\|_{L^1(B_R)} &\leq C\left( R^{n(1-\frac{1}{1+\varepsilon})+1} \|D_x^2 u\|_{L^{1+\varepsilon}(B_R)} + R^{-1} \|u\|_{L^1(B_R)}\right) \\
							%&\leq C\left( R^{n(1-\frac{1}{1+\varepsilon})+1} \|\Delta u\|_{L^{1+\varepsilon}(B_R)} + R^{-1}R^{n-\alpha}\right)  \\
							&\leq C\left( R^{n(1-\frac{1}{1+\varepsilon})+1} R^{\frac{n-p\beta}{1+\varepsilon}} +R^{-1}R^{n-\alpha}\right)  \\
							&\leq C R^{n+1-\frac{\alpha+2}{1+\varepsilon}},
\end{align*}
The second estimate in \eqref{3rdEstimate} can be obtained by a similar process. Last, the fact that $n-(p+1)\beta<0$ gives
\begin{align*}
	\|D_x^2 u\|_{L^k(B_R)}^k &\leq C \left( \int_{B_{2R}} |\Delta u|^kdx + R^{n-k(n+2)}(\int_{B_{2R}} |u|dx)^k\right) \\
			&\leq C\left( \int_{B_{2R}} v^{p+1}dx + R^{n-k(n+2)} R^{k(n-\alpha)}\right) \\
			&\leq C\left( F(2R) + R^{n-(p+1)\beta}\right) \\
			&\leq C F(2R),
\end{align*}
and hence the first estimate in \eqref{4thEstimate} follows, and similarly we get the second estimate.
$\Box$

\textbf{Proof of Proposition \ref{estimateOnSphere}:} By Lemma \ref{sphereEstimate}, $\exists \tilde{R}\in[R,2R]$, Proposition \ref{estimateOnSphere} follows from Proposition \ref{estimateOnBR} immediately. $\Box$

\begin{lemma}\label{F4R}
There exists $M>0$ such that $\exists \{R_j\}\rightarrow +\infty$,
\begin{align*}
	F(4R_j) \leq M F(R_j).
\end{align*}
\end{lemma}
Proof. Suppose not, then for any $M>0$ and any $\{R_j\}\rightarrow+\infty$, we have
\begin{align*}
	F(4R_j) > M F(R_j).
\end{align*}
Take $M>5^n$ and $R_{j+1}= 4R_j$ with $R_0>1$. Therefore,
\begin{align*}
	F(R_j) > M^j F(R_0),
\end{align*}
which leads to a contradiction with $F(R_j)\leq CR_j^n \leq C(4^jR_0)^n$. $\Box$

\section{Proof of Liouville Theorem}

From now on, without loss of generality, we may assume $p\geq q.$ By Lemma \ref{comparisonPrinciple}, $\|v\|^{p+1}_{L^{p+1}(B_R)}\leq \|u\|^{q+1}_{L^{q+1}(B_R)}$. By the Rellich-Poho\v{z}aev type identity in Lemma \ref{pohozaevId}, we can denote
\begin{align}\label{FR}
    F(R):=\int_{B_R}u^{q+1}\leq CG_1(R)+CG_2(R),
\end{align}
where
\begin{align}
   & G_1(R)=R^n\int_{\mathbb S^{n-1}}u^{q+1}(R)d\sigma, \label{g1}  \\
    & G_2(R)=R^n\int_{\mathbb S^{n-1}}(|D_x u(R)|+R^{-1}u(R))(|D_x v(R)|+R^{-1}v(R))d\sigma.  \label{g2}
\end{align}
Heuristically, we are aiming for estimate as %a sequence of $\{R_j\}\rightarrow+\infty$ such that
\begin{align}\label{goalEstimate}
    %G_i(R_j)\leq CR_j^{-a_i}F^{1-\delta_i}(R_j), \quad i=1,2,
	G_i(R)\leq CR^{-a_i}F^{1-\delta_i}(4R), \quad i=1,2.
\end{align}
%and $a_i>0$ and $\delta_i>0$.
Then by Lemma \ref{F4R} there exists a sequence $\{R_j\}\rightarrow+\infty$ such that
%$F(R_j)\leq CR_j^{-\bar{a}}F^{\bar{b}}(4R_j)\leq CR_j^{-\bar{a}}F^{\bar{b}}(R_j)$.
\begin{align*}
	G_i(R_j)\leq CR^{-a_i}F^{1-\delta_i}(R_j), \quad i=1,2.
\end{align*}
Suppose there are infinitely many $R_j$'s such that $G_1(R_j)\geq G_2(R_j)$, then take that subsequence of $\{R_j\}$ and still denote as $\{R_j\}$. We do the same if there are infinitely many $R_j$'s such that $G_1(R_j)\leq G_2(R_j)$. So, there are only two cases we shall deal with: there exists a sequence $\{R_j\}\rightarrow+\infty$ such that
\begin{enumerate}
	\item[Case] 1: $G_1(R_j)\geq G_2(R_j)$. Then we prove $a_1>0$, $\delta_1>0$. So, $F^{\delta_1}(R_j)\leq CR_j^{-a_1}\rightarrow 0$,
	\item[Case] 2: $G_1(R_j)\leq G_2(R_j)$. Then we prove $a_2>0$, $\delta_2>0$. So, $F^{\delta_2}(R_j)\leq CR_j^{-a_2}\rightarrow 0$.
\end{enumerate}
Then we conclude that $F(R)\equiv 0$.

\begin{comment}
Let $\bar{a}=\min\{a,\tilde{a}\}$ and $\bar{b}=\max\{b,\tilde{b}\}$. Suppose we have
\begin{align}\label{aBarbBar}
	\bar{a}>0, \ \bar{b}<1,
\end{align}
and we show that
\begin{align*}
	F(R)\leq F(\tilde{R}) \leq CR^{-\bar{a}}F^{\bar{b}}(4R), \quad R\geq 1,
\end{align*}
implies $F(R)\equiv 0$.

First, we claim that there exists $M>0$ such that $\exists \{R_j\}$ with $R_j\rightarrow +\infty$,
\begin{align*}
	F(4R_j) \leq M F(R_j).
\end{align*}
Suppose not, then for any $M>0$ and any $\{R_j\}\rightarrow+\infty$ we have
\begin{align*}
	F(4R_j) > M F(R_j).
\end{align*}
Take $M>5^n$ and $R_{j+1}= 4R_j$ with $R_0>1$. Therefore,
\begin{align*}
	F(R_j) > M^j F(R_0),
\end{align*}
which leads to a contradiction with $F(R_j)\leq CR_j^n \leq C(4^jR_0)^n$.
\end{comment}

Surprisingly, for both cases $a_i\approx(\alpha+\beta+2-n)\delta_i$. Indeed, we have
\begin{theorem}\label{liouvilleThm2}
For $F(R)$ defined as \eqref{FR} and $\alpha,\beta$ defined as \eqref{scalingComponent}, there exists a sequence $\{R_j\}\rightarrow+\infty$ such that
\begin{align*}
	F(R_j) \leq C R_j^{-(\alpha+\beta+2-n)+o(1)}.
\end{align*}
\end{theorem}
Hence, Theorem \ref{liouvilleThm} is a direct consequence of Theorem \ref{liouvilleThm2}, and we only need to prove Theorem \ref{liouvilleThm2} for case 1 and 2.

\subsection{Case 1: Estimate for $G_1(R)$}

%\subsection{$G_1$ and $G_2$}
%By Lemma \ref{comparisonPrinciple} and the assumption in Theorem \ref{liouvilleThm} we have
%\begin{align}\label{vEstimate}
	%\int_{B_R} v^s \leq CR^{n-s\beta}.
%\end{align}

According to previous discussion in the introduction, we assume that
\begin{align*}
	p\geq q >0, \quad pq>1, \quad \beta\leq \alpha<n-2, \quad
	n\geq 3,
\end{align*}
hence in particular
\begin{align}\label{pLowerBound}
	p>\frac{n}{n-2}.
\end{align}

\begin{remark}
For systems \eqref{laneEmdenSystem2} with double bounded coefficients, \eqref{pLowerBound} is a necessary %s and sufficient
condition for existence of positive solution, see \cite{LL13}.
\end{remark}

%In addition, we assume that $0\leq n-s\beta<1$ since if $n-s\beta < 0$, Liouville type of theorems can be obtained easily. %For the case $n-s\beta=0$, the solutions have finite energy, and then they have to be radial and therefore must be 0, cf. \cite{ChenLi2009}.
%Also, we may assume $s \geq p$. If $s<p$, \eqref{vEstimateBetter} can be interpolated from energy estimates in Lemma \ref{energyEstimates} by H\"{o}lder inequality. So, let $\tilde{s}=\max\{p,s\}$, and we replace $s$ in \eqref{vEstimateBetter} with $\tilde{s}$.
In addition to our assumption that $n-s\beta<1$,
since we have energy inequalities \eqref{uvEsitmatePQ}, we can assume $s\geq p$.
Also, if $n-s\beta < 0$, \eqref{vEstimateBetter} implies $v\equiv 0$ and hence $u\equiv 0$. So, we assume $n-s\beta \geq 0$.
Let $l=\frac{s}{p}$, then
\begin{align}\label{rangeL}
	l\geq 1, \ \text{and } \frac{n-1}{p\beta} < l \leq \frac{n}{p\beta}.
\end{align}

%We are going to prove Theorem \ref{liouvilleThm2} for this $s$ (or $l$). In fact,
It is worthy to point out that, what the proof of Lane-Emden conjecture really needs is a ``breakthrough" on the energy estimate \eqref{uvEsitmatePQ}. $s$ in \eqref{vEstimateBetter} needs not be very large but enough to satisfy $n-s\beta<1$. In other words, $s$ can be very close to $\frac{n-1}{\beta}$, and it is sufficient to prove Theorem \ref{liouvilleThm}. %Also, notice that, Lane-Emden conjecture in the case $n=3,4$ follows from Theorem \ref{liouvilleThm} directly by letting $s=p$.

The strategies of attacking $G_1$ and $G_2$ are the same. Basically, first by H\"{o}lder inequality we split the quantities on sphere $\mathbb S^{n-1}$ into two parts. One has a lower (than original) index after embedding, and the other has a higher one. Then we estimate the latter part by $F(R)$, and thus we get a feedback estimate as \eqref{goalEstimate}.
%{\bf{Step 2.}} Estimate for $G_2(R)$.

Let
 \begin{align*}
 	k=\frac{p+1}{p}.
 \end{align*}
  Since $p\beta=\alpha+2$, $n-(p+1)\beta=n-2-(\alpha+\beta)<0$ by \eqref{subcriticalRegion2}. Thus, $n-kp\beta<0$ as $n-lp\beta\geq 0$, it follows that $l<k$.
   %Let $\lambda$ and $\mu$ be positive numbers such that
  % \begin{align}\label{embeddingIndex}
%       \ \text{and }
  %      \frac{1}{\mu} \geq \frac{1}{k} - \frac{2}{n-1}.
   %\end{align}

 {\bf{Subcase 1.1}} $\frac{1}{l}  \geq \frac{2}{n-1} + \frac{1}{q+1}$.

Note that in this subcase, since $l\geq 1$, we must have $n\geq 4$ (i.e., $n\neq 3$).
By \eqref{pLowerBound} we see that $k=1+\frac 1 p<1+\frac{n-2}{n} = \frac 2 n (n-1)\leq \frac{n-1}{2}$. Take
\begin{align*}
	\frac{1}{\mu} = \frac{1}{k} - \frac{2}{n-1}.
\end{align*}
So, $W^{2,k}(\mathbb S^{n-1})\hookrightarrow L^{\mu}(\mathbb S^{n-1})$.

Take
\begin{align*}
	\frac{1}{\lambda} = \frac{1}{l} - \frac{2}{n-1} \geq \frac{1}{q+1}.
\end{align*}
Then $W^{2,l+\varepsilon}(\mathbb S^{n-1})\hookrightarrow L^{\lambda}(\mathbb S^{n-1})$.

Direct verification shows that $ \frac{1}{\mu}= \frac{1}{k} - \frac{2}{n-1} \leq \frac{1}{q+1}$ which is due to \eqref{subcriticalRegion},
so we have %$\frac{1}{\mu}\leq \frac{1}{q+1}\leq \frac{1}{\lambda}$.
 \begin{align*}
 	\frac{1}{\mu}\leq \frac{1}{q+1}\leq \frac{1}{\lambda}.
 \end{align*}
 Then by H\"{o}lder inequality and Sobolev embedding \eqref{sobolevEmbedding}, we have (with notation \eqref{sphereNorm})
 \begin{align}
     \|u\|_{q+1}(R) &\leq \|u\|_{\lambda}^{\theta}\|u\|_{\mu}^{1-\theta}(R) \label{uHolderInequality} \\
                 &\leq C(R^2\|D_x^2 u\|_{l+\varepsilon}(R)+\|u\|_1(R))^{\theta}(R^2\|D_x^2 u\|_k(R)+\|u\|_1(R))^{1-\theta}, \label{uInequality}
 \end{align}
 where $\theta\in [0,1]$ and
 \begin{align} \label{theta}
     \frac{1}{q+1}=\frac{\theta}{\lambda}+\frac{1-\theta}{\mu}.
 \end{align}

   %then we can solve $\theta$ by \eqref{embeddingIndex}
  %\begin{align*}
  %	\theta = \frac{\frac{1}{1+q}-\frac{1}{k}+\frac{2}{n-1}}{\frac{1}{l}-\frac{1}{k}},
  %\end{align*}
  Since $l$ can be 1 (then $W^{2,p}$-estimate fails for $\|D_x^2 u\|_1(R)$), we add an $\varepsilon$ to $l$ for later use of $W^{2,p}$-estimate. $\varepsilon$ can be any real positive number and later will be chosen sufficiently small.

  To get desired estimate, we have requirements in form of inequalities involving parameters, such as $\alpha,\beta,\varepsilon$ and etc. To verify those requirements very often we just verify strict inequalities with $\varepsilon=0$ because such inequalities continuously depend on $\varepsilon$. % See section 3.2.

 So, by \eqref{g1} and \eqref{uInequality}
  \begin{align}\label{g1Estimate1}
  	G_1(R) \leq CR^n \left( (R^2\|D_x^2 u\|_{l+\varepsilon}(R)+\|u\|_1(R))^{\theta}(R^2\|D_x^2 u\|_k(R)+\|u\|_1(R))^{1-\theta}\right) ^{q+1}.
  \end{align}
Then by Proposition \ref{estimateOnSphere}, there exists $\tilde{R}\in[R,2R]$ such that
\begin{align*}
    G_1(\tilde{R})
    &\leq CR^n\left(( R^2R^{-\frac{lp\beta}{l+\varepsilon}}+R^{-2-\alpha})^{\theta}(R^2(R^{-n}F(4R))^{\frac{1}{k}}  +R^{-\alpha})^{1-\theta} \right) ^{q+1}\\
    &\leq CR^n\left( R^{2-\frac{lp\beta\theta}{l+\varepsilon}-\frac{n(1-\theta)}{k}}F^{\frac{1-\theta}{k}}(4R)\right) ^{q+1} \\
    &\leq R^{-a_1}F^{1-\delta_1}(4R),
\end{align*}
where the last inequality is due to $R^{-\frac{n}{k}}>R^{-\alpha-2}$ and
\begin{align}\label{abG1}
    a_1 &= a_1^{\varepsilon}=(q+1)(\frac{lp\beta\theta}{l+\varepsilon}+\frac{np(1-\theta)}{p+1}-2-\frac{n}{1+q}),\\
    1-\delta_1 &= \frac{(1-\theta)p(q+1)}{p+1}.
\end{align}
 % Since for sufficiently small $\varepsilon$, $a_{\varepsilon}, \tilde{a}_{\varepsilon}>0$ and $b,\tilde{b}<1$ are just a perturbation of
 Since for sufficiently small $\varepsilon$, $a_1^{\varepsilon}>0$ and $\delta_1>0$ are just a perturbation of
 \begin{align}\label{a1delta1}
 %	a_0,\tilde{a}_{0}>0, \ \text{and } b,\tilde{b}<1,
 	a_1^0>0, \ \text{and } \delta_1>0,
 \end{align}
we only need to prove \eqref{a1delta1} is true.

  Since $lp=s$, $p\beta=\alpha+2$ and $q\alpha=\beta+2$,
  \begin{align*}
  	a_1^0 &= p\beta\theta(q+1) + (1-\delta_1)n - 2(q+1)-n \\
  	   %&= p\beta\theta(q+1) - 2(q+1)+(b-1)n \\
  	   &= (q+1)(p\beta\theta-2) -\delta_1n \\
  	   &= (q+1)(p\beta(\theta-1)+p\beta -2) -\delta_1n \\
  	   %&= (q+1)(-\frac{2}{pq-1}b(p+1)+\alpha) +(b-1)n\\
  	   &= (q+1)(-\alpha (1-\delta_1)+\alpha)-\delta_1n \\
  	   &= \delta_1((q+1)\alpha-n) \\
  	   &= (\alpha + \beta +2 -n)\delta_1.
  \end{align*}
So we just need to prove $\delta_1>0$. By \eqref{abG1} and \eqref{theta} we have
 \begin{align*}
 	% \Leftrightarrow&
 	  &(1-\theta)p(q+1) < p+1 \\
 	 \Leftrightarrow& \frac{\frac{1}{l}-\frac{2}{n-1}-\frac{1}{1+q}}{\frac{1}{l}-\frac{1}{k}}(q+1)<k \\
 	 \Leftrightarrow& (\frac{1}{l}-\frac{2}{n-1})(q+1) -1< \frac{k}{l}-1 \\
 	 \Leftrightarrow& \frac{1}{l}(q+1-1-\frac{1}{p}) < \frac{2}{n-1}(q+1) \\
 	 \Leftrightarrow& \frac{pq-1}{s} < \frac{2(q+1)}{n-1} \\
 	 \Leftrightarrow& n-1 < s\beta,
 \end{align*}
 and the last inequality is included in our assumption. So, we have proved subcase 1.1.

{\bf{Subcase 1.2}} $\frac{1}{l}  < \frac{2}{n-1}  + \frac{1}{q+1}$.

As discussed in the beginning of subcase 1.1, $k<\frac{n-1}{2}$ if $n>3$. Since $l<k$, $\frac{1}{l}  > \frac{2}{n-1}$ for $n>3$. When $n=3$, since $l\geq 1$ by \eqref{rangeL},  $\frac{1}{l} \leq 1 =\frac{2}{n-1}$. %Therefore, $\frac{1}{l}  \geq \frac{2}{n-1}$ for $n\geq 3$.

Therefore, for $n>3$, take
\begin{align*}
	\frac{1}{\lambda} = \frac{1}{l} - \frac{2}{n-1} < \frac{1}{q+1}, %\ \text{if } n>3,
\end{align*}
 and for $n=3$, take
 \begin{align*}
 	\lambda = \infty, % \ \text{if } n=3,
\end{align*}
so we have
\begin{align*}
	W^{2,l+\varepsilon}(\mathbb S^{n-1})\hookrightarrow L^{\lambda}(\mathbb S^{n-1}), \ n\geq 3.
\end{align*}

 So, % let $\theta = 1$, and \eqref{G1Inequality} still holds, i.e.
 \begin{align*}
 	\|u\|_{q+1}(R) \leq C\|u\|_{\lambda}(R)\leq C(R^2\|D_x^2 u\|_{l+\varepsilon}(R)+\|u\|_1(R) ).
 \end{align*}
 Therefore, by Proposition \ref{estimateOnSphere} there exists $\tilde{R}\in[R,2R]$ such that
 \begin{align}
   	G_1(\tilde{R}) &\leq CR^n  (R^2\|D_x^2 u\|_{l+\varepsilon}(R)+\|u\|_1(R))^{q+1} \\
   							&\leq CR^{n}(R^2 R^{-\frac{lp\beta}{l+\varepsilon}}+R^{-\alpha})^{q+1}  \\
   							&\leq CR^{n+(2-\frac{lp\beta}{l+\varepsilon})(q+1) }.
   \end{align}
 So,
 \begin{align*}
 	F(\tilde{R}) &\leq C R^{n+(2-\frac{lp\beta}{l+\varepsilon})(q+1) } \\
 			&\leq CR^{n+(2-p\beta)(q+1)+\frac{\varepsilon p\beta}{l+\varepsilon}(q+1)} \\
 			&\leq CR^{-(\alpha+\beta+2-n) + \frac{\varepsilon p\beta}{l+\varepsilon}(q+1)}.
 \end{align*}
 Since $\varepsilon$ can be arbitrarily small,
 \begin{align*}
 	F(\tilde{R}) \leq CR^{-(\alpha+\beta+2-n) + o(1)}.
 \end{align*}

Thus, we have proved Case 1.

\subsection{Case 2: Estimate for $G_2(R)$.}

Let
\begin{align*}
	 m=\frac{q+1}{q}.
\end{align*}

%As discussed in Step 1, $\frac{1}{l}  \geq \frac{2}{n-1}$ for $n\geq 3$.

{\bf{Subcase 2.1}} $m<n-1$. %, i.e. $q>\frac{1}{n-2}$.

With $z,z'>0$ and $\frac 1 z + \frac{1}{z'}=1$, \eqref{g2} becomes,
\begin{align}\label{G2holderInequality}
\left.\begin{array}{ll}
	G_2(R) &\leq CR^n\||D_x u| + R^{-1}u\|_z\||D_x v| + R^{-1}v\|_{z'}(R) \\
				&\leq CR^n(\|D_x u\|_z(R)+R^{-1}\|u\|_z(R))(\|D_x v\|_{z'}(R) + R^{-1}\|v\|_z(R)) \\
				&\leq CR^n(\|D_x u\|_z(R)+R^{-1}\|u\|_1(R))(\|D_x v\|_{z'}(R) + R^{-1}\|v\|_1(R)),
\end{array}
 \right.
\end{align}
where the last inequality is due to
\begin{align*}
	\|u\|_z(R) \leq C(R\|D_x u\|_z(R) + \|u\|_1(R)), \ \text{and } \|v\|_{z'}(R) \leq C(R\|D_x v\|_{z'}(R) + \|v\|_1(R)).
\end{align*}
%For sufficiently small $\varepsilon$ which will be chosen later, and
Assume there exists $z$ (we shall check the existence later) such that by Sobolev Embedding \eqref{sobolevEmbedding},
\begin{align}
&\left.\begin{array}{ll}\label{g2UEmbedding}
		\|D_x u\|_z(R) &\leq \|D_x u\|_{\rho_1}^{\tau_1}(R)\|D_x u\|_{\gamma_1}^{1-\tau_1}(R) \\
							&\leq C(R\|D^2_x u\|_{l+\varepsilon}(R) + \|D_x u\|_1(R))^{\tau_1}(R\|D^2_x u\|_k(R)+ \|D_x u\|_1(R))^{1-\tau_1},
\end{array}
 \right.  \\
& \left.\begin{array}{ll}\label{g2VEmbedding}
 	\|D_x v\|_{z'}(R) &\leq \|D_x v\|_{\rho_2}^{\tau_2}(R)\|D_x v\|_{\gamma_2}^{1-\tau_2}(R) \\
 							&\leq C(R\|D^2_x v\|_{1+\varepsilon}(R) + \|D_x v\|_1(R))^{\tau_2}(R\|D^2_x v\|_m(R)+ \|D_x v\|_1(R))^{1-\tau_2},		
 \end{array}
  \right. 			
\end{align}
where $\tau_1,\tau_2\in[0,1]$ and
\begin{align}
  	\frac{1}{z} = \frac{\tau_1}{\rho_1} + \frac{1-\tau_1}{\gamma_1}, \label{zDef} \\
  	\frac{1}{z'} = \frac{\tau_2}{\rho_2} + \frac{1-\tau_2}{\gamma_2},	\label{zDefPrime}
\end{align}
and since $l<k\leq m<n-1$, %plus if $q>\frac{1}{n-2}$ (i.e. $\frac{1}{m}>\frac{1}{n-1}$) %when $n>3$ and $l=1$ when $n=3$,
define
\begin{align}
	\frac{1}{\rho_1} &= \frac{1}{l} - \frac{1}{n-1}, \ \ \frac{1}{\gamma_1} = \frac{1}{k} - \frac{1}{n-1}, \\
		\frac{1}{\rho_2} &= 1 - \frac{1}{n-1}, \ \ \frac{1}{\gamma_2} = \frac{1}{m} - \frac{1}{n-1}. \label{rhoGammaDef}
\end{align}
So, we have
\begin{align*}
W^{1,l+\varepsilon}(\mathbb{S}^{n-1})\hookrightarrow L^{\rho_1}(\mathbb{S}^{n-1}), \ W^{1,k}(\mathbb{S}^{n-1})\hookrightarrow L^{\gamma_1}(\mathbb{S}^{n-1}), \\
W^{1,1+\varepsilon}(\mathbb{S}^{n-1})\hookrightarrow L^{\rho_2}(\mathbb{S}^{n-1}), \ W^{1,m}(\mathbb{S}^{n-1})\hookrightarrow L^{\gamma_2}(\mathbb{S}^{n-1}).
\end{align*}
To verify the existence of such $z$, by \eqref{zDef}-\eqref{rhoGammaDef}, we expect that
\begin{align}\label{rangeZ}
	\max\left\lbrace \frac{1}{k} - \frac{1}{n-1}, \frac{1}{n-1}\right\rbrace \leq \frac{1}{z} \leq \min \left\lbrace \frac{1}{l}-\frac{1}{n-1}, \frac{1}{q+1}+\frac{1}{n-1}\right\rbrace.
\end{align}
Thus, we need to verify, (i) $\frac{1}{k} - \frac{1}{n-1}\leq\frac{1}{l}-\frac{1}{n-1}$, (ii) $\frac{1}{n-1}\leq\frac{1}{l}-\frac{1}{n-1}$, (iii) $\frac{1}{n-1}\leq\frac{1}{q+1}+\frac{1}{n-1}$, (iv) $\frac{1}{k}-\frac{1}{n-1}\leq\frac{1}{q+1}+\frac{1}{n-1}$.

Since $l<k$, (i) is true. (ii) holds for $n> 3$ as discussed at the beginning of subcase 1.2 $\frac 1 l>\frac 1 k>\frac{2}{n-1}$; for $n=3$, take $s=p$ and then $l=1$, so (ii) still holds. (iii) is obvious. (iv) is equivalent to $\frac{1}{p+1}+\frac{1}{q+1}\geq 1-\frac{2}{n-1}$, which is guaranteed by \eqref{subcriticalRegion}.
\begin{comment}
\begin{align*}
	\frac{1}{k} - \frac{1}{n-1}<\frac{1}{l}-\frac{1}{n-1}.
\end{align*}
Meanwhile,
\begin{align*}
	\frac{1}{k}-\frac{1}{n-1}<\frac{1}{q+1}+\frac{1}{n-1}
	\Leftrightarrow \frac{1}{p+1}+\frac{1}{q+1}> 1-\frac{2}{n-1},
\end{align*}
which is guaranteed by \eqref{subcriticalRegion}.
Last,
\begin{align*}
	\frac{1}{n-1}  \leq \frac{1}{l} - \frac{1}{n-1}
\end{align*}
holds for $n\geq 3$ as discussed at the beginning of Step 2.
\end{comment}

So, we put \eqref{g2UEmbedding} and \eqref{g2VEmbedding} in \eqref{G2holderInequality} and get
\begin{align}\label{g2Estimate}
\left. \begin{array}{ll}
 G_2(R)\leq & C R^{n+2}(\|D_x^2 u\|_{l+\varepsilon}(R) + R^{-1}\|D_x u\|_1(R) + R^{-2}\|u\|_1(R))^{\tau_1} \\
             & \times(\|D_x^2 u\|_k(R) + R^{-1}\|D_x u\|_1(R) + R^{-2}\|u\|_1(R))^{1-\tau_1} \\
             & \times(\|D_x^2 v\|_{1+\varepsilon}(R) + R^{-1}\|D_x v\|_1(R) + R^{-2}\|v\|_1(R))^{\tau_2} \\
             & \times(\|D_x^2 v\|_m(R) + R^{-1}\|D_x v\|_1(R) + R^{-2}\|v\|_1(R))^{1-\tau_2}.
\end{array}\right.
\end{align}

Then by Proposition \ref{estimateOnSphere}, there exists $\tilde{R}\in[R,2R]$ such that
\begin{align*}
	G_2(\tilde{R}) &\leq C R^{n+2} R^\frac{-p\beta\tau_1}{1+\varepsilon/l} \left( (R^{-n}F(4R))^{\frac{1}{k}}+  R^{-\frac{\alpha+2}{1+\varepsilon}} +R^{-\alpha-2}\right) ^{1-\tau_1}  \\
	&\times R^{-\frac{q\alpha\tau_2}{1+\varepsilon/l}} \left( (R^{-n}F(4R))^{\frac{1}{m}} +R^{-\frac{\beta+2}{1+\varepsilon}}+R^{-\beta-2}\right)^{1-\tau_2} \\
	&\leq C R^{-a_2^{\varepsilon}}F^{1-\delta_2}(4R),
\end{align*}
where the last inequality is due to $R^{-\frac{n}{k}}>R^{-\alpha-2}$ and $R^{-\frac{n}{m}}>R^{-\beta-2}$. Meanwhile,
\begin{align}
	a_2&= a_2^{\varepsilon} = -n-2 +\frac{p\beta\tau_1}{1+\varepsilon/l} + \frac{q\alpha\tau_2}{1+\varepsilon/l} + n\frac{1-\tau_1}{k} + n\frac{1-\tau_2}{m}, \\
	1-\delta_2 &= \frac{1-\tau_1}{k} + \frac{1-\tau_2}{m}. \label{delta2}
\end{align}
Similar to subcase 1.1, we only need to prove
\begin{align*}
	a_2^0 >0, \ \delta_2>0.
\end{align*}
Surprisingly, similar to $a_1\approx (\alpha+\beta+2-n)\delta_1$, we have $a_2\approx (\alpha+\beta+2-n)\delta_2$ since we can prove $a_2^0 =(\alpha+\beta+2-n)\delta_2$. Indeed,
\begin{align*}
	a^0_2 &= -n-2 +p\beta(\tau_1-1)+p\beta+q\alpha(\tau_2-1)+q\alpha+n(1-\delta_2) \\
					&= -n-2 - p\beta k \frac{1-\tau_1}{k} - q\alpha m\frac{1-\tau_2}{m} + \alpha+\beta+4 +n(1-\delta_2) \\
					&= \alpha + \beta +2-n - (\alpha+\beta+2) (1-\delta_2) + n(1-\delta_2)\\
					&= (\alpha+\beta+2-n)\delta_2,
\end{align*}
where the third equality above is due to $p\beta k=(p+1)\beta=(q+1)\alpha=q\alpha m$ and $(p+1)\beta=\alpha+\beta+2$.
So, we only need to prove $\delta_2 >0$ or equivalently by \eqref{zDef}, \eqref{zDefPrime} and \eqref{delta2},
\begin{align}
(m-\frac{k}{l})\frac{1}{z}+(\frac{k}{n-1}+(m-1)(k-1))\frac{1}{l}+\frac{m-2}{n-1}-(m-1)>0, \label{zAndl}
\end{align}
To achieve this, we take the upper bound of $\frac{1}{z}$ in \eqref{rangeZ} and see whether \eqref{zAndl} holds.

{\bf{Case 2.1.1}} If $\frac{1}{l}-\frac{1}{n-1} \geq \frac{1}{q+1}+\frac{1}{n-1}$, then let $\frac{1}{z}=\frac{1}{q+1}+\frac{1}{n-1}$, and \eqref{zAndl} becomes,
\begin{align*}
	&(\frac{1}{pq}-\frac{p+1}{p(q+1)})\frac{1}{l}+\frac{q+1}{q}(\frac{1}{n-1}+\frac{1}{q+1})+\frac{1-q}{(n-1)q}-\frac{1}{q}>0 \\
	&\Leftrightarrow (\frac{1}{pq}-\frac{p+1}{p(q+1)})\frac{1}{l}+ \frac{2}{q(n-1)} > 0 \\
	&\Leftrightarrow -\frac{2}{\beta s} + \frac{2}{n-1} > 0 \\
	&\Leftrightarrow s\beta > n-1.
\end{align*}
{\bf{Case 2.1.2}} If $\frac{1}{l}-\frac{1}{n-1} < \frac{1}{q+1}+\frac{1}{n-1}$, then let $\frac{1}{z}=\frac{1}{l}-\frac{1}{n-1}$, and \eqref{zAndl} becomes, %A trick in calculation here is that let $l=l_0=\frac{n-1}{p\beta}$, if the inequality holds, then for $l$ near $l_0$ it still holds.
\begin{align*}
	&(m-\frac{k}{l})(\frac{1}{l}-\frac{1}{n-1})+(\frac{k}{n-1}+(m-1)(k-1))\frac{1}{l}+\frac{m-2}{n-1}-(m-1)>0 \\
  &\Leftrightarrow -\frac{k}{l^2}+(m+\frac{k}{n-1}+\frac{k}{n-1}+(m-1)(k-1))\frac{1}{l}+\frac{m-2}{n-1}-(m-1)>0 \\
  %&\Leftrightarrow -\frac{k}{l^2}+(\frac{q+1}{q}+2\frac{p+1}{(n-1)p}+\frac{1}{pq})\frac{1}{l}>\frac{2}{n-1}+\frac{1}{q} \\
 % &\Leftrightarrow -\frac{k}{l^2}+(\frac{p(q+1)}{(p+1)q}+\frac{2}{n-1}+\frac{1}{(p+1)q})\frac{p+1}{lp}>\frac{2}{n-1}+\frac{1}{q} \\
  &\Leftrightarrow -\frac{k}{l^2}+(\frac{p}{p+1}+\frac{1}{q}+\frac{2}{n-1})\frac{k}{l}>\frac{2}{n-1}+\frac{1}{q} \\
%\end{align*}
%Since $k=\frac{p+1}{p}$, then the inequality above becomes
%\begin{align*}
	&\Leftrightarrow - \frac{k}{l^2} + (1+k(\frac{2}{n-1}+\frac{1}{q}))\frac{1}{l} - (\frac{2}{n-1}+\frac{1}{q}) >0 \\
	&\Leftrightarrow (\frac{k}{l} -1)(\frac{1}{l} -(\frac{2}{n-1}+\frac{1}{q})) < 0 \\
	&\Leftrightarrow \frac{1}{k} < \frac{1}{l} < \frac{2}{n-1}+\frac{1}{q}.
\end{align*}
Notice that $\frac{1}{l} < \frac{2}{n-1}+\frac{1}{q}$ holds under the assumption of case 2.1.2, and $\frac{1}{k} < \frac{1}{l}$ since $l<k$. %Also, it is easy to verify $\frac{1}{k} <\frac{2}{n-1}+\frac{1}{q+1}<\frac{2}{n-1}+\frac{1}{q}$.
In all, \eqref{zAndl} always holds under our assumption $n-s\beta<1$.

{\bf{Subcase 2.2}} $m\geq n-1$. %, i.e. $q \leq \frac{1}{n-2}$.

First, we have for any $\gamma\in[1,\infty)$,
\begin{align*}
W^{1,m}(\mathbb{S}^{n-1}) \hookrightarrow L^{\gamma}(\mathbb{S}^{n-1}).
\end{align*}
Then we claim $\frac{1}{l}>\frac{1}{n-1}$. Suppose $\frac{1}{l}\leq \frac{1}{n-1}$, then $k>l\geq n-1$, hence $p \leq \frac{1}{n-2}$, which is not possible due to \eqref{pLowerBound}.
Take $\frac{1}{z}=\frac{1}{l}-\frac{1}{n-1}$ %and any $\gamma_2\in[z',+\infty)$, and we have
then
\begin{align*}
W^{1,l+\varepsilon}(\mathbb{S}^{n-1}) \hookrightarrow L^{z}(\mathbb{S}^{n-1}).
\end{align*}
Therefore, by Sobolev embedding and \eqref{G2holderInequality}
\begin{align*}
	G_2(R) &\leq  CR^n(\|D_x u\|_z(R)+R^{-1}\|u\|_1(R))(\|D_x v\|_{z'}(R) + R^{-1}\|v\|_1(R)) \\
			&\leq C R^{n+2}(\|D_x^2 u\|_{l+\varepsilon} + R^{-1}\|D_x u\|_1 + R^{-2}\|u\|_1)(\|D_x^2 v\|_m + R^{-1}\|D_x v\|_1 + R^{-2}\|v\|_1).
\end{align*}
Similarly to previous work, there exists a $\tilde{R}\in[R,2R]$ such that % \eqref{g2Estimate} becomes,
\begin{align*}
	G_2(\tilde{R}) &\leq C R^{n+2} R^\frac{-p\beta}{1+\varepsilon/l}
 \left( (R^{-n}F(4R))^{\frac{1}{m}} +R^{-\frac{\beta+2}{1+\varepsilon}}+R^{-\beta-2}\right) \\
	&\leq C R^{-a_2^{\varepsilon}}F^{1-\delta_2}(4R),
\end{align*}
where
\begin{align}\label{abG2}
	a_2 &= a_2^{\varepsilon} = -n-2 +\frac{p\beta}{1+\varepsilon/l} + \frac{n}{m}, \\
	1-\delta_2 &= \frac{1}{m}.
\end{align}
Direct verification shows that
\begin{align*}
	a_2^0 = (\alpha+\beta+2-n)\delta_2,
\end{align*}
and obviously $\delta_2>0$ so $\alpha_2^0>0$.

Thus, we have proved Case 2.

%\subsection{Estimate on $\mathbb S^{n-1}$}
%Here are some observations:
 %\item By \eqref{holderIndex} we see that there is only one free variable to choose since $k$ is fixed;
  %\item $\theta$ increases as $\lambda$ increases. %, and $\lambda$ increases as $l$ increases.
%\end{itemize}

%Thus, the estimate on $G_1, G_2$ relies on all those quantities on $\mathbb S^{n-1}$.

\newpage

\providecommand{\bysame}{\leavevmode\hbox to3em{\hrulefill}\thinspace}
\providecommand{\MR}{\relax\ifhmode\unskip\space\fi MR }
\providecommand{\MRhref}[2]{%
  \href{http://www.ams.org/mathscinet-getitem?mr=#1}{#2}
}
\providecommand{\href}[2]{#2}

\noindent\href{mailto:ze.cheng@colorado.edu}{ze.cheng@colorado.edu}\\
\href{mailto:genggenghuang1986@gmail.com}{genggenghuang1986@gmail.com}\\
\href{mailto:congmingli@gmail.com}{congmingli@gmail.com}
\end{document}